\documentclass{amsart}
\usepackage{amssymb}
\usepackage{pstricks,pst-node}

\renewcommand{\Im}{\operatorname{Im}}

\newcommand{\SL}{\operatorname{SL}}

\newcommand{\Res}{\operatorname{Res}}

\newcommand{\Ga}{\Gamma}

\newcommand{\al}{\alpha}
\newcommand{\be}{\beta}
\newcommand{\ga}{\gamma}
\newcommand{\de}{\delta}
\newcommand{\ze}{\zeta}
\newcommand{\eps}{\varepsilon}

\newcommand{\ka}{\kappa}
\newcommand{\La}{\Lambda}
\newcommand{\la}{\lambda}
\newcommand{\sig}{\sigma}
\renewcommand{\phi}{\varphi}
\newcommand{\om}{\omega}

\newcommand{\C}{\mathbf C}
\newcommand{\N}{\mathbf N}
\newcommand{\Z}{\mathbf Z}

\renewcommand{\P}{\mathbb P}

\newcommand{\A}{{\mathcal A}}
\newcommand{\B}{{\mathcal B}}
\newcommand{\CC}{{\mathcal C}}

\renewcommand{\L}{{\mathcal L}}
\newcommand{\M}{{\mathcal M}}

\renewcommand{\o}{\otimes}

\newcommand{\inv}{^{-1}}

\newcommand{\pmat}[1]{\begin{pmatrix}#1\end{pmatrix}}
\newcommand{\vmat}[1]{\begin{vmatrix}#1\end{vmatrix}}

\newcommand{\exch}{\leftrightarrow}

\newcommand{\To}{\Rightarrow}
\newcommand{\from}{\leftarrow}
\newcommand{\From}{\Leftarrow}

\newcommand{\spmat}[1]{\left(\begin{smallmatrix}#1\end{smallmatrix}\right)}

\theoremstyle{definition}
\newtheorem{df}{Definition}[section]
\newtheorem{rem}[df]{Remark}
\newtheorem{infrem}[df]{Informal remark}

\theoremstyle{plain}
\newtheorem{theo}{Theorem}

\newtheorem{prop}[df]{Proposition}

\numberwithin{equation}{section}

\title{Quantum $\SL(3,\C)$'s: the missing case}
\author{Christian Ohn}
\date{October 15, 2002}
\address{Universit\'e de Valenciennes\\ Laboratoire de Math\'ematiques,
I.S.T.V.\\Le Mont Houy\\F-59313 Valenciennes Cedex 9\\
France}
\email{christian.ohn@univ-valenciennes.fr}
\thanks{2000 \emph{Mathematics Subject Classication.} Primary 20G42;
Secondary 14M15, 16S37, 16S38.}

\begin{document}

\begin{abstract}
We study the only missing case in the classification of quantum $\SL(3,\C)$'s
undertaken in our paper [\emph{J.\ of Algebra} {\bf213} (1999), 721--756],
thereby completing this classification.
\end{abstract}

\maketitle

\section*{Introduction}

The aim of this paper is to complete the classification of quantum
$\SL(3)$'s undertaken in \cite{qsl3}.

Roughly speaking, we call a quantum $\SL(3)$ any Hopf algebra (over $\C$)
whose (finite-dimensional) comodules are ``similar'' to the modules of the
(ordinary) group $\SL(3,\C)$ (see Definition~\ref{Gqdef}). Given such a Hopf
algebra $\A$, it will have, among other things, two simple (nonisomorphic)
comodules $V,W$ of dimension~$3$, and the usual decomposition rules for
tensor products will imply the existence of $\A$-comodule morphisms
\eqref{basic}, satisfying certain compatibility conditions \eqref{coh}. The
tuple $\L_\A$, consisting of those two comodules and eight morphisms, will be
called the \emph{basic quantum $\SL(3)$ datum} (BQD for short) associated to
$\A$.

Conversely, starting from an ``abstract'' BQD $\L$ (i.e.\ where $V,W$ are
just vector spaces and \eqref{basic} just linear maps satisfying the
aforementioned conditions), we may reconstruct a Hopf algebra $\A_\L$ via the
usual Tannakian procedure.

We may now state our main result.
\begin{theo}
\label{maintheo}
~
\begin{enumerate}
\item If $\A$ is a quantum $\SL(3)$, then $\L_\A$ is a BQD.
\item If $\L$ is a BQD, then $\A_\L$ is a quantum $\SL(3)$.
\item The correspondences $\A\mapsto\L_\A$ and $\L\mapsto\A_\L$ are inverse
of each other between quantum $\SL(3)$'s (up to Hopf algebra isomorphism) and
BQDs (up to equivalence; see Definition~\ref{BQDdef}).
\item BQDs can be explicitly classified up to equivalence, yielding a
classification of quantum $\SL(3)$'s up to Hopf algebra isomorphism.
\end{enumerate}
\end{theo}

An almost complete version of this theorem was stated and proved in
\cite{qsl3}. More precisely, we found a classification of all BQDs, except
for one class related to elliptic curves, and we proved the theorem for all
BQDs outside this class. In the present paper, we settle the study of this
last class, thereby proving Theorem~A in full.

Moreover, a geometric analysis of this class of BQDs yields the following
contribution to a question raised in the Introduction of \cite{AST}.
\begin{theo}
\label{ASTtheo}
Let $\B$ be the quantum three-space associated to an arbitrary quantum
$\SL(3)$. Then $\B$ cannot be a Sklyanin algebra. In other words, the scheme
of point modules \cite{ATV,ATV2} of $\B$ (which is a cubic divisor in $\P^2$)
cannot be an elliptic curve.
\end{theo}

The paper is organized as follows. After some recollections from \cite{qsl3}
in Section~\ref{recollsec} (definition of a quantum $\SL(3)$ and of a BQD,
and the correspondences between them), Section~\ref{ellsec} recalls and
studies the form of the only class of BQDs not covered by the results of
\cite{qsl3}; this will finish the classification of BQDs. In
Section~\ref{shapesec}, we introduce the shape algebra \cite[Section~5]{qsl3}
for this class of BQDs and we determine the associated flag variety (in the
sense of \cite{qflag}). These geometric data suggest to view the shape
algebra as a twist (in the sense of \cite{Zh}) of another algebra, which we
show in Section~\ref{twistshapesec} to be isomorphic to the shape algebra of
another BQD, already covered by the results of \cite{qsl3}. In
Section~\ref{proofmainsec}, we finish the proof of Theorem~\ref{maintheo} by
carrying over the necessary properties from the untwisted to the twisted
shape algebra. Theorem~\ref{ASTtheo} is proved in Section~\ref{proofASTsec}
by picking up some leftovers from Section~\ref{ellsec}.

In the Appendix, we show a result on twists of Koszul algebras that is needed
in Section~\ref{proofmainsec}, but may also be of independent interest.

{\bf Conventions.}
We denote by $\Z$ (resp.~$\N$, $\C$) the set of integers (resp.\ nonnegative
integers, complex numbers). All vector spaces, algebras and tensor products
are over $\C$.

\section{Recollections from \cite{qsl3}}
\label{recollsec}

Recall that the group $\SL(3)$ is linearly reductive and that its simple
modules are parametrized by their highest weights, which are pairs
$(k,l)\in\N^2$. Recall further that the dimension of the simple module of
highest weight $(k,l)$ is given by $d_{(k,l)}:=(k+1)(l+1)(k+l+2)/2$. For any
$\la,\mu,\nu\in\N^2$, denote by $m_{\la\mu}^\nu$ the multiplicity of the
simple module of highest weight $\nu$ inside the tensor product of those of
(respective) highest weights $\la$ and $\mu$. (Recall also that
$m_{\la\mu}^\nu$ can, in principle, be determined in a purely combinatorial
way.)

Our main objects of interest may now be defined as follows (see \cite{qsl3}).
\begin{df}
\label{Gqdef}
We call a \emph{quantum $\SL(3)$} any (not necessarily commutative) Hopf
algebra $\A$ (over $\C$) such that
\begin{enumerate}
\item there is a family $\{V_\la\mid\la\in\N^2\}$ of simple and pairwise
nonisomorphic $\A$-comodules, with $\dim V_\la=d_\la$,
\item every $\A$-comodule is isomorphic to a direct sum of these,
\item for every $\la,\mu\in\N^2$, $V_\la\o V_\mu$ is isomorphic to
$\bigoplus_{\nu} m_{\la\mu}^\nu V_\nu$.
\end{enumerate}
\end{df}
In particular, if we write $V:=V_{(1,0)}$ and $W:=V_{(0,1)}$, then
Condition~(c) implies the existence of $\A$-comodule
morphisms 
\begin{equation}
\label{basic}
\begin{aligned}
A&:V\o V\to W,
&\qquad a&:W\to V\o V
\\
B&:W\o W\to V,
&b&:V\to W\o W
\\
C&:W\o V\to\C,
&c&:\C\to V\o W
\\
D&:V\o W\to\C,
&d&:\C\to W\o V,
\end{aligned}
\end{equation}
each being unique up to a scalar. We showed in \cite[Propositions 3.1~and
3.2]{qsl3} that these maps must, for an appropriate choice of these scalars,
satisfy the following compatibility conditions:
\begin{subequations}
\label{coh}
\begin{align}
(1_V\o C)(c\o1_V)&=1_V,&(D\o1_V)(1_V\o d)&=1_V\\
Aa&=1_W,&&\\
C(A\o1_V)&=\om\,D(1_V\o A),&(1_V\o a)c&=\om\,(a\o 1_V)d\\
\om\,(C\o 1_V)(1_W\o a)&=B,&(A\o 1_W)(1_V\o c)&=b\\
Dc&=\kappa\,1_\C,&Cd&=\kappa\,1_\C
\end{align}
\begin{align}
(1_V\o A)(a\o 1_V)(A\o 1_V)(1_V\o a)&=\rho\,(1_{V\o W}+cD)\\
(A\o 1_V)(1_V\o a)(1_V\o A)(a_V\o 1)&=\rho\,(1_{W\o V}+dC),
\end{align}
\end{subequations}
where
\begin{itemize}
\item $\om$ is a 3rd root of unity,
\item $\ka=q^{-2}+1+q^2$ and $\rho=(q+q\inv)^{-2}$ for some $q\in\C$, $q\ne0$,
with $q^2$ either $1$ or not a root of unity.
\end{itemize}
\begin{df}
\label{BQDdef}
A \emph{basic quantum $\SL(3)$ datum} (BQD for short) is a tuple
$\L=(V,W,A,a,B,b,C,c,D,d)$ consisting of two vector spaces $V$ and $W$ of
dimension~$3$, and of eight linear maps \eqref{basic}, satisfying
Conditions~\eqref{BQDdef} (with $\om,\ka,q$ as indicated).
Two BQDs are called \emph{equivalent} if one can be obtained from the other
through (any combination of) base change, rescaling of the maps
\eqref{basic}, and interchanging $V\exch W$, $A\exch B$, $a\exch b$, $C\exch
D$, $c\exch d$.
\end{df}
Thus, each quantum $\SL(3)$, $\A$, gives rise to a BQD $\L_\A$ (which is
really defined only up to equivalence).

Conversely, start with a BQD $\L$, and define an algebra $\A_\L$ with
$(9+9)$ generators $t^i_j$ ($i,j=1,2,3$) and $u^\al_\be$ ($\al,\be=1,2,3$),
and relations
\begin{align*}
A_{ij}^\al\,t^i_kt^j_\ell&=u^\al_\be\,A^\be_{k\ell},
&
t^i_kt^j_\ell\,a^{k\ell}_\be&=a^{ij}_\al\,u^\al_\be
\\
B_{\al\be}^i\,u^\al_\ga u^\be_\de&=t^i_j\,B^j_{\ga\de},
&
u^\al_\ga u^\be_\de\,b^{\ga\de}_j&=b^{\al\be}_i\,t^i_j
\\
C_{\al i}\,u^\al_\be t^i_j&=C_{\be j},
&
t^i_ju^\al_\be\,c^{j\be}&=c^{i\al}
\\
D_{i\al}\,t^i_ju^\al_\be&=D_{j\be},
&
u^\al_\be t^i_j\,d^{\be j}&=d^{\al i}.
\end{align*}
(Here, we have chosen bases $x_1,x_2,x_3$ of
$V$ and $y_1,y_2,y_3$ of $W$, and set $A(x_i\o x_j)=A_{ij}^\al y_\al$,
etc.) Then $\A_\L$ possesses a Hopf algebra structure given by
\begin{equation}
\label{costr}
\begin{aligned}
\Delta(t^i_j)&=t^i_k\o t^k_j,
&\Delta(u^\al_\be)&=u^\al_\ga\o u^\ga_\be
\\
\eps(t^i_j)&=\de^i_j,
&\eps(u^\al_\be)&=\de^\al_\be
\\
S(t^i_j)&=c^{i\be}\,u^\al_\be\,C_{\al j},
&S(u^\al_\be)&=d^{\al j}\,t^i_j\,D_{i\be}.
\end{aligned}
\end{equation}
The main difficulty is to prove that the Hopf algebra $\A_\L$ is indeed a
quantum $\SL(3)$ in the sense of Definition~\ref{Gqdef}. We were able to do
this in \cite{qsl3} for all BQDs, except for one class, which we will
describe in the next section.

\section{Classification of Case I.h}
\label{ellsec}

If we set $Q^i_j:=c^{i\al}D_{j\al}$, then (\ref{coh}a) implies
$(Q^{-1})^i_j=d^{\al i}C_{\al_j}$. Now it follows from \eqref{costr} that
\[
S^2(t^i_j)=Q^i_kt^k_l(Q^{-1})^l_j,
\]
so $S^2$ is ``encoded'' by the linear map $Q:V\to V$.

In \cite[Section~10]{qsl3}, we used the possible Jordan normal forms of $Q$
and the value of $\om$ as first criteria for the classification of BQDs. One
possibility, called \emph{Type I} in \cite{qsl3}, consists in taking $Q=1_V$
and $\om=1$. The condition $Q=1_V$ amounts to setting $W=V^*$, with $C,c,D,d$
the obvious canonical maps. (In particular, we now have $S^2=1_\A$,
reflecting the fact that $W$ is both the left dual and the right dual of
$V$.) In this case, (\ref{coh}a) and (\ref{coh}e) are automatically satisfied
(for $\ka=3$, so $q^2=1$ and $\rho=\frac14$).

Now consider the ``quantum determinants,'' appearing in (\ref{coh}c):
\begin{align*} e:=(1_V\o a)c=\om(a\o1_V)d&:\C\to V\o V\o V \\
E:=C(A\o1_V)=\om D(1_V\o A)&:V\o V\o V\to\C. \end{align*} If $\L$ is of
Type~I, then Conditions~(\ref{coh}c) imply that \[ e=\la+s,\qquad E=\La+S, \]
with $\la,\La$ totally antisymmetric and $s,S$ totally symmetric. In
particular, choosing dual bases $x_i$ in $V$ and $y_\al$ in $W(=V^*)$, we may
view $s$ and $S$ as homogeneous polynomials of degree~$3$, i.e.\ as cubic
curves in the projective plane $\P^2$ and in its dual plane ${\P^2}^*$,
respectively (unless $s=0$ or $S=0$).

Using the standard classification of cubic curves in $\P^2$, we were then
able in \cite{qsl3} to classify all possible forms of $e$ and $E$ that
satisfy Conditions~(\ref{coh}bfg), except in the case where $s$ is an
elliptic curve.

In this particular case, called \emph{Case I.h} in \cite{qsl3}, the bases of
$V$~and $W$ may be chosen in such a way that $a$ reads
\begin{align*}
a:\;&y_1\mapsto\al\,x_2\o x_3+\be\,x_3\o x_2+\ga\,x_1\o x_1
\\
&y_2\mapsto\al\,x_3\o x_1+\be\,x_1\o x_3+\ga\,x_2\o x_2
\\
&y_3\mapsto\al\,x_1\o x_2+\be\,x_2\o x_1+\ga\,x_3\o x_3,
\end{align*}
with $\ga\ne0$ and $(\al+\be)^3+\ga^3\ne0$ (so that $s$ is indeed elliptic).
It then follows from (\ref{coh}b) that $A$ must read
\begin{align*}
A:\;x_1\o x_1&\mapsto\ga'y_1,&
x_1\o x_2&\mapsto\al'y_3,&
x_1\o x_3&\mapsto\be'y_2
\\
x_2\o x_1&\mapsto\be'y_3,&
x_2\o x_2&\mapsto\ga'y_2,&
x_2\o x_3&\mapsto\al'y_1
\\
x_3\o x_1&\mapsto\al'y_2,&
x_3\o x_2&\mapsto\be'y_1,&
x_3\o x_3&\mapsto\ga'y_3,
\end{align*}
with
\begin{equation}
\label{sp}
\al\al'+\be\be'+\ga\ga'=1.
\end{equation}
Moreover, Conditions~(\ref{coh}fg) now read
\begin{equation}
\label{p2p2}
\begin{aligned}
P_0&:=\al^2{\al'}^2+\be^2{\be'}^2+\ga^2{\ga'}^2
-2\,\al\al'\be\be'-2\,\al\al'\ga\ga'-2\,\be\be'\ga\ga'=0
\\
P_1&:=\al^2\be'\ga'+\be^2\al'\ga'+\ga^2\al'\be'=0
\\
P_2&:={\al'}^2\be\ga+{\be'}^2\al\ga+{\ga'}^2\al\be=0.
\end{aligned}
\end{equation}

\begin{prop}
\label{elldegen}
Up to base change, we may assume that $\al=\al'=0$.
\end{prop}

\begin{proof}
\emph{Case 1} : $\al=0$ or $\al'=0$. Substituting into $P_2$ or
$P_1$ shows that $\al=\al'=0$.

\emph{Case 2} : $\be=0$ or $\be'=0$. Permuting two of the three basis vectors
takes us back to Case~1.

\emph{Case 3} : $\ga'=0$ (recall that $\ga=0$ has been ruled out by
ellipticity of $s$). Substituting into $P_1$ shows that $\al'=0$ or $\be'=0$,
which takes us back to Case~1 or Case~2.

\emph{Case 4} : none of $\al,\al',\be,\be',\ga,\ga'$ equals zero. Viewing
$P_0,P_1,P_2$ as polynomials in $\al'$, their resultants must vanish:
\begin{align*}
0&=\Res_{\al'}(P_0,P_1)
=\vmat{\al^2&-2\al(\be\be'+\ga\ga')&(\be\be'-\ga\ga')^2
\\
\be^2\ga'+\ga^2\be'&\al^2\be'\ga'&0
\\
0&\be^2\ga'+\ga^2\be'&\al^2\be'\ga'}
\\
&=\be^2\ga^4{\be'}^4+2(\al^3+\be^3-\ga^3)\be\ga^2{\be'}^3\ga'
\\
&\qquad
+(\al^6+\be^6+\ga^6+2\al^3\be^3+2\al^3\ga^3-4\be^3\ga^3){\be'}^2{\ga'}^2
\\
&\qquad+2(\al^3+\ga^3-\be^3)\be^2\ga\be'{\ga'}^3+\be^4\ga^2{\ga'}^4
\\
&=:Q_1
\\
0&=\Res_{\al'}(P_1,P_2)
=\vmat{\be^2\ga'+\ga^2\be'&\al^2\be'\ga'&0
\\
0&\be^2\ga'+\ga^2\be'&\al^2\be'\ga'
\\
\be\ga&0&\al({\be'}^2\ga+{\ga'}^2\be)}
\\
&=\al\ga^5{\be'}^4+2\al\be^2\ga^3{\be'}^3\ga'
+(\al^3+\be^3+\ga^3)\al\be\ga{\be'}^2{\ga'}^2
+2\al\be^3\ga^2\be'{\ga'}^3+\al\be^5{\ga'}^4
\\
&=:Q_2.
\end{align*}
(These resultants make sense, because the leading coefficients of
$P_0,P_1,P_2$ are nonzero: in particular, if we had $\be^2\ga'+\ga^2\be'=0$,
then substituting into $P_1$ would imply $\al=0$, $\be'=0$, or $\ga'=0$.)

Viewing $Q_1,Q_2$ as polynomials in $\be'$, their resultant must again vanish:
\begin{align}
0&=\Res_{\be'}(Q_1,Q_2)=\text{(a $6\times6$ determinant)}\notag
\\
\label{ellnotell}
&={\ga'}^{16}(\al\be\ga)^{10}
\Bigl[(\al^3+\be^3+\ga^3)^3-(3\al\be\ga)^3\Bigr].
\end{align}
(The author confesses not to have computed this $6\times6$ determinant
by hand!)

Therefore,
$\al^3+\be^3+\ga^3=3\ze\al\be\ga$ for some 3rd root of unity $\ze$, i.e.
\[
(\ze\ga+\al+\be)(\ze\ga+j\al+j^2\be)(\ze\ga+j^2\al+j\be)=0,
\]
where $j$ is a primitive 3rd root of unity. But $\ze\ga+\al+\be=0$ is ruled
out by ellipticity of the curve $s$, so we may assume that
$\ze\ga+j\al+j^2\be=0$ (exchanging $j\exch j^2$ if necessary).

Now consider the following change of basis in $V$:
\[
\left\{
\begin{aligned}
x'_1&=\ze x_1+x_2+x_3
\\
x'_2&=\ze x_1+jx_2+j^2x_3
\\
x'_3&=\ze x_1+j^2x_2+jx_3
\end{aligned}
\right.
\]
(together with the dual basis $y'_1,y'_2,y'_3$ in $W$). Then the map $a$ reads
\begin{align*}
(3\ze)\;a:\;&y'_1\mapsto(\ze\ga+j^2\al+j\be)\,x'_3\o
x'_2+(\ze\ga+\al+\be)\,x'_1\o x'_1
\\
&y'_2\mapsto(\ze\ga+j^2\al+j\be)\,x'_1\o x'_3+(\ze\ga+\al+\be)\,x'_2\o x'_2
\\
&y'_3\mapsto(\ze\ga+j^2\al+j\be)\,x'_2\o x'_1+(\ze\ga+\al+\be)\,x'_3\o x'_3,
\end{align*}
so we are taken back to Case~1.
\end{proof}

Substituting $\al=\al'=0$ back into \eqref{sp} and \eqref{p2p2} now yields
$\be\be'=\ga\ga'=\frac12$. Now we still have one degree of freedom to rescale
the maps $a$ and $A$, so we are left with one essential parameter, say
\[
t:=-\frac\ga\be=-\frac{\be'}{\ga'}.
\]
This completes the classification of BQDs, undertaken in
\cite[Section~10]{qsl3}.

\section{The shape algebra and its flag variety}
\label{shapesec}

Given a BQD $\L$, recall \cite[Section~5]{qsl3} that its \emph{shape algebra}
$\M_\L$ is generated by $x_i$ ($i=1,2,3$) and $y_\al$ ($\al=1,2,3$), with
defining relations
\begin{align*}
a^{ij}_\al x_ix_j&=0,&c^{i\al}x_iy_\al&=0
\\
b^{\al\be}_iy_\al y_\be&=0,&y_\al x_i&=-(q+q^{-1})a^{jk}_\al
A^\be_{ki}x_jy_\be
\end{align*}
This algebra has a natural $\N^2$-grading
$\M_\L=\bigoplus_{(k,l)\in\N^2}V_{(k,l)}$, with the $x_i$ of degree $(1,0)$
and the $y_\al$ of degree $(0,1)$. (Moreover, it carries a natural
$\A_\L$-comodule algebra structure, with each $V_{(k,l)}$ becoming an
$\A_\L$-comodule: these are the natural candidates to show that $\A_\L$
satisfies Definition~\ref{Gqdef}.)

When $\L$ is the BQD of Case~I.h (with $\al=\al'=0$ and
$\be\be'=\ga\ga'=\frac12$; cf.\ Section~\ref{ellsec}), the defining relations
of $\M_\L$ read
\begin{equation}
\label{ellshape}
\begin{gathered}
\begin{aligned}
x_3x_2&=t\,x_1^2,
\\
x_1x_3&=t\,x_2^2,
\\
x_2x_1&=t\,x_3^2,
\end{aligned}
\qquad
\begin{aligned}
t\,y_2y_3&=y_1^2
\\
t\,y_3y_1&=y_2^2
\\
t\,y_1y_2&=y_3^2
\end{aligned}
\\[2mm]
\begin{aligned}
y_1x_1&=x_2y_2,
&y_2x_1&=t\,x_2y_3,
&t\,y_3x_1&=x_2y_1
\\
t\,y_1x_2&=x_3y_2,
&y_2x_2&=x_3y_3,
&y_3x_2&=t\,x_3y_1
\\
y_1x_3&=t\,x_1y_2,
&t\,y_2x_3&=x_1y_3,
&y_3x_3&=x_1y_1
\end{aligned}
\\[2mm]
x_1y_1+x_2y_2+x_3y_3=0.
\end{gathered}
\end{equation}
Now let us determine the flag variety of $\M_\L$ as defined in \cite{qflag}.
To do this, modify Relations~\eqref{ellshape} as follows: in each
relation, adorn the right factor of each term with a $'$ if the left factor
is an $x_i$, and with a $''$ if the left factor is a $y_\al$ (e.g.\ the
first relation on the fourth line now reads $y_1x''_1=x_2y'_2$).

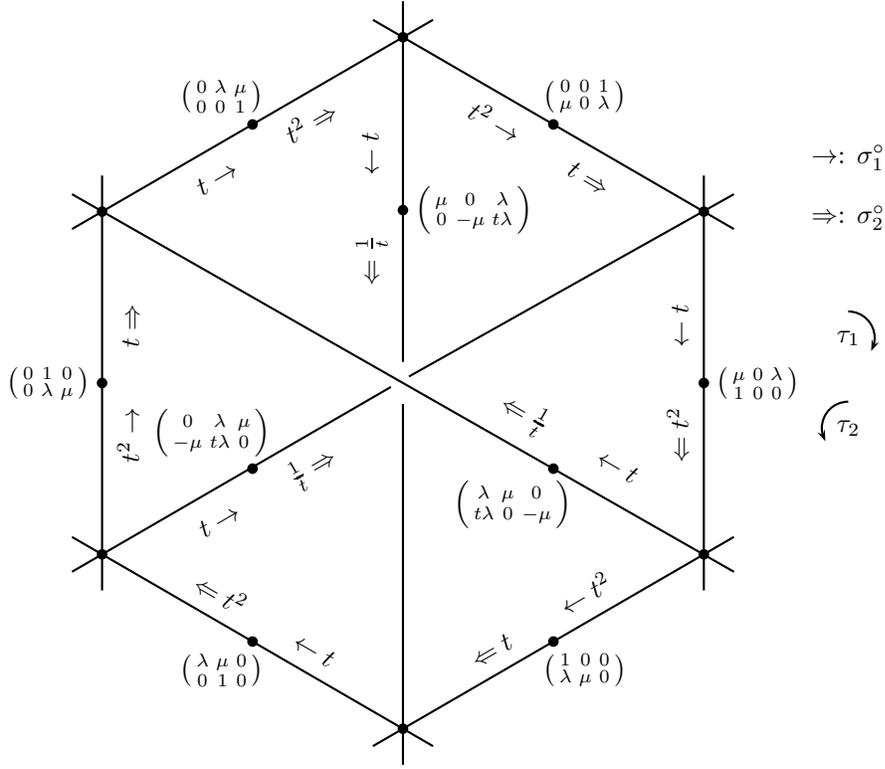
\begin{figure}
\begin{center}
\psset{unit=.8mm}
\begin{pspicture}(-5,-6)(105,120)
\pscircle*(50,0){2pt}
\pscircle*(100,29){2pt}
\pscircle*(100,86){2pt}
\pscircle*(50,115){2pt}
\pscircle*(0,86){2pt}
\pscircle*(0,29){2pt}
\qline(45,-2.9)(105,31.9)
\qline(100,23)(100,92)
\qline(105,83.1)(45,117.9)
\qline(55,117.9)(-5,83.1)
\qline(0,92)(0,23)
\qline(-5,31.9)(55,-2.9)
\qline(-5,26.1)(48,56.84)
\qline(51,58.58)(105,88.9)
\qline(-5,88.9)(105,26.1)
\qline(50,-6)(50,54)
\qline(50,61)(50,121)
\pscircle*(25,14.5){2pt}
\uput{2}[240]{0}(25,14.5){$\spmat{\la&\mu&0\\0&1&0}$}
\pscircle*(75,14.5){2pt}
\uput{2}[300]{0}(75,14.5){$\spmat{1&0&0\\{}\la&\mu&0}$}
\pscircle*(100,57.5){2pt}
\uput{2}[0]{0}(100,57.5){$\spmat{\mu&0&\la\\1&0&0}$}
\pscircle*(75,100.5){2pt}
\uput{2}[60]{0}(75,100.5){$\spmat{0&0&1\\{}\mu&0&\la}$}
\pscircle*(25,100.5){2pt}
\uput{2}[120]{0}(25,100.5){$\spmat{0&\la&\mu\\0&0&1}$}
\pscircle*(0,57.5){2pt}
\uput{2}[180]{0}(0,57.5){$\spmat{0&1&0\\0&\la&\mu}$}
\pscircle*(75,43.25){2pt}
\uput{2}[240]{0}(75,43.25){$\spmat{\la&\mu&0\\t\la&0&-\mu}$}
\pscircle*(50,86.25){2pt}
\uput{2}[0]{0}(50,86.25){$\spmat{\mu&0&\la\\0&-\mu&t\la}$}
\pscircle*(25,43.25){2pt}
\uput{2}[120]{0}(25,43.25){$\spmat{0&\la&\mu\\-\mu&t\la&0}$}
\pcline[linestyle=none](50,0)(100,29)
\aput{:U}{$\From t\qquad\from t^2$}
\pcline[linestyle=none](100,26)(100,89)
\aput{:U}{$\From t^2\qquad\from t$}
\pcline[linestyle=none](50,115)(100,86)
\bput{:U}{$t^2\to\qquad t\To$}
\pcline[linestyle=none](0,86)(50,115)
\bput{:U}{$t\to\qquad t^2\To$}
\pcline[linestyle=none](0,29)(0,86)
\bput{:U}{$t^2\to\qquad t\To$}
\pcline[linestyle=none](0,29)(50,0)
\aput{:U}{$\From t^2\qquad\from t$}
\pcline[linestyle=none](0,29)(50,58)
\bput{:U}{$t\to\qquad\frac1t\To$}
\pcline[linestyle=none](50,58)(100,29)
\aput{:U}{$\From\frac1t\qquad\from t$}
\pcline[linestyle=none](50,58)(50,115)
\aput{:U}{$\From\frac1t\qquad\from t$}
\psarc{<-}(124,65){10pt}{-30}{90}
\psarc{->}(124,50){10pt}{90}{210}
\rput{0}(124,65){$\tau_1$}
\rput{0}(124,50){$\tau_2$}
\rput{0}(124,95){$\to$: $\sig_1^\circ$}
\rput{0}(124,85){$\To$: $\sig_2^\circ$}
\end{pspicture}
\end{center}
\caption{Flag variety for Case I.h, with a generic point on each irreducible
component (top row in each matrix: coordinates in $\P^2$; bottom row:
coordinates in ${\P^2}^*$). Single and double arrows represent $\sig_1^\circ$
and $\sig_2^\circ$, respectively; each arrow multiplies $\mu$ by the
coefficient next to it. The automorphisms $\tau_1,\tau_2$ are ``rotations by
120~degrees'' in the directions shown.}
\label{Ihpicture}
\end{figure}

View $(x_1:x_2:x_3)$ and $(y_1:y_2:y_3)$ as homogeneous coordinates in
$\P^2$ and in ${\P^2}^*$, respectively. Since Relations~\eqref{ellshape}
are $\N^2$-homogeneous, their modified version may now be seen as
defining equations for a subscheme
\[
\Gamma\subset(\P^2\times{\P^2}^*)\times(\P^2\times{\P^2}^*)
\times(\P^2\times{\P^2}^*).
\]
One easily checks that $\Gamma$ is of the form
\[
\Gamma=\{(p,\sig_1(p),\sig_2(p))\mid p\in X\}
\]
for some subscheme $X\subset\P^2\times{\P^2}^*$ and some automorphisms
$\sig_1,\sig_2$ of $X$: indeed, the scheme $X$ is the variety with nine
irreducible components pictured in Figure~\ref{Ihpicture},
and the automorphisms $\sig_1,\sig_2$ naturally decompose into
$\sig_1=\tau_1\sig_1^\circ$, $\sig_2=\tau_2\sig_2^\circ$, where 
$\tau_1,\tau_2$ are the automorphisms of $\P^2\times{\P^2}^*$ given by
\[
\tau_1:\pmat{x_1&x_2&x_3\\y_1&y_2&y_3}\mapsto\pmat{x_3&x_1&x_2\\y_3&y_1&y_2},
\qquad
\tau_2:\pmat{x_1&x_2&x_3\\y_1&y_2&y_3}\mapsto\pmat{x_2&x_3&x_1\\y_2&y_3&y_1},
\]
and where $\sig_1^\circ,\sig_2^\circ$ are also described in
Figure~\ref{Ihpicture}, e.g.\
$\sig_1^\circ:\spmat{0&1&0\\0&\la&\mu}\mapsto\spmat{0&1&0\\0&\la&t^2\mu}$ and
$\sig_2^\circ:\spmat{0&1&0\\0&\la&\mu}\mapsto\spmat{0&1&0\\0&\la&t\mu}$.
\begin{infrem}
Since $\M_\L$ is a quantum analogue of the multihomogeneous coordinate ring
of the (ordinary) flag variety of $\SL(3)$, nontrivial characters of $\M_\L$
should, in some sense, correspond to ``quantum Borel subgroups'' of
$\A_\L$. On the other hand, such characters obviously correspond to
simultaneous fixed points of $\sig_1,\sig_2$. But in the case considered
here, there is no such fixed point: so once Theorem~\ref{maintheo} will be
proved, we will have a quantum group with the same representations as
$\SL(3)$, although in some sense, it has \emph{no} quantum Borel subgroup to
induce them from!
\end{infrem}

\section{Twisting the shape algebra}
\label{twistshapesec}

The natural decompositions $\sig_i=\tau_i\sig_i^\circ$ suggest to consider
$\M_\L$ as the twist $R^\tau$ (in the sense of \cite{Zh}) of another
$\N^2$-graded algebra $R$, whose associated flag variety will be
$(X,\sig_1^\circ,\sig_2^\circ)$ instead of $(X,\sig_1,\sig_2)$ (this works by
an obvious multigraded version of \cite[Proposition~8.9]{ATV2}, noting that
the $\tau_1,\tau_2$ commute with each other and with
$\sig_1^\circ,\sig_2^\circ$). The algebra $R$ will also be generated by $x_i$
($i=1,2,3$) and $y_\al$ ($\al=1,2,3$), with defining relations obtained by
twisting \eqref{ellshape} backwards:
\begin{equation}
\label{ellshapetwist}
\begin{gathered}
\begin{aligned}
x_3x_1&=t\,x_1x_3,
\\
x_1x_2&=t\,x_2x_1,
\\
x_2x_3&=t\,x_3x_2,
\end{aligned}
\qquad
\begin{aligned}
t\,y_2y_1&=y_1y_2
\\
t\,y_3y_2&=y_2y_3
\\
t\,y_1y_3&=y_3y_1
\end{aligned}
\\[2mm]
\begin{aligned}
y_1x_2&=x_2y_1,
&y_2x_2&=t\,x_2y_2,
&t\,y_3x_2&=x_2y_3
\\
t\,y_1x_3&=x_3y_1,
&y_2x_3&=x_3y_2,
&y_3x_3&=t\,x_3y_3
\\
y_1x_1&=t\,x_1y_1,
&t\,y_2x_1&=x_1y_2,
&y_3x_1&=x_1y_3
\end{aligned}
\\[2mm]
x_1y_3+x_2y_1+x_3y_2=0.
\end{gathered}
\end{equation}
A straightforward computation shows that relations \eqref{ellshapetwist},
with $y_1,y_2,y_3$ renamed to $y_2,y_3,y_1$, turn out to define the shape
algebra $\M_{\L^\circ}$, where $\L^\circ$ is the BQD of Type~I defined by the
maps
\begin{gather*}
\begin{aligned}
a:\;&y_1\mapsto\la\,x_2\o x_3+\mu\,x_3\o x_2
\\
&y_2\mapsto\la\,x_3\o x_1+\mu\,x_1\o x_3
\\
&y_3\mapsto\la\,x_1\o x_2+\mu\,x_2\o x_1
\end{aligned}
\\[2mm]
\begin{aligned}
A:\;x_1\o x_1&\mapsto0,&
x_1\o x_2&\mapsto\la'y_3,&
x_1\o x_3&\mapsto\mu'y_2
\\
x_2\o x_1&\mapsto\mu'y_3,&
x_2\o x_2&\mapsto0,&
x_2\o x_3&\mapsto\la'y_1
\\
x_3\o x_1&\mapsto\la'y_2,&
x_3\o x_2&\mapsto\mu'y_1,&
x_3\o x_3&\mapsto0
\end{aligned}
\end{gather*}
(Case I.e in \cite[Section~10]{qsl3}), with $\la\la'=\mu\mu'=\frac12$ and
$t=-\frac\mu\la=-\frac{\la'}{\mu'}$.

\section{Proof of Theorem~\ref{maintheo}}
\label{proofmainsec}

In \cite{qsl3}, all intermediate results in the (almost complete) proof of
Theorem~\ref{maintheo} are valid for an arbitrary BQD $\L$, except for
\cite[Proposition~5.3]{qsl3} (and the resulting \cite[Corollary~5.4]{qsl3}),
which relies on a case by case analysis that is not valid in Case~I.h.
More precisely, the proof of Theorem~\ref{maintheo}
will be complete if we show the following facts for this particular case:
\begin{enumerate}
\item the shape algebra $\M_\L$ is a Koszul algebra (when viewed as an
$\N$-graded algebra via the total grading),
\item $\dim V_{(k,l)}=d_{(k,l)}$ for all $(k,l)\in\N^2$.
\end{enumerate}
Since the BQD $\L^\circ$ introduced in Section~\ref{twistshapesec} is already
covered by \cite{qsl3}, Properties (a) and (b) are true for $\M_{\L^\circ}$.

Therefore, $\M_\L$, being a twist of $\M_{\L^\circ}$, also satisfies those
two properties: for Property~(b), this is automatic, and for Property~(a), it
follows from Proposition~\ref{Koszultwist} in the Appendix.
Theorem~\ref{maintheo} follows.

\section{Proof of Theorem~\ref{ASTtheo}}
\label{proofASTsec}

To each BQD $\L$ are associated two ``quantum three-spaces,'' namely the
quadratic algebras $\B_\L:=T(V)/(\Im a)$ and $\CC_\L:=T(W)/(\Im b)$. For
Case~I.h, $\B_\L$ is defined by the following relations:
\begin{align*}
\al\,x_2x_3+\be\,x_3x_2+\ga\,x_1^2&=0
\\
\al\,x_3x_1+\be\,x_1x_3+\ga\,x_2^2&=0
\\
\al\,x_1x_2+\be\,x_2x_1+\ga\,x_3^2&=0,
\end{align*}
(The defining relations of $\CC_\L$ are similar.) This algebra is one of the
regular algebras of dimension~$3$ studied in \cite{AS}, where the matrix $Q$
(see Section~\ref{ellsec}) also played a role in the classification of such
algebras. Algebras for which $Q=1$ (called of \emph{Type~A} in \cite{AS})
give rise to the cubic curve $s$ in $\P^2$ as in Section~\ref{ellsec}, so let
us call this curve the \emph{AS-curve}. The AS-curve associated to $\B_\L$ is
given by
\[
\ga(x_1^3+x_2^3+x_3^3)=3(\al+\be)x_1x_2x_3.
\]
Another cubic curve in $\P^2$ has been associated to $\B_\L$ in \cite{ATV},
namely the scheme of its point modules, given by
\[
(\al\be\ga)(x_1^3+x_2^3+x_3^3)=(\al^3+\be^3+\ga^3)x_1x_2x_3.
\]
Call this curve the \emph{ATV-curve} associated to $\B_\L$.

\begin{rem}
The ATV-curve is defined for an arbitrary regular algebra of dimension~$3$,
whereas the AS-curve is defined only for those of Type~A. Note however, as
observed in \cite{ATV}, that even when both curves are defined, they \emph{do
not} coincide in general!
\end{rem}

Recall \cite{ATV,OF} that $\B_\L$ is called a \emph{Sklyanin algebra} if its
ATV-curve is elliptic. But by Proposition~\ref{elldegen}, this is impossible:
the ATV-curve of $\B_\L$ must degenerate to a triangle (which may also be
viewed as the image of the flag variety $X$ under the projection
$\P^2\times{\P^2}^*\to\P^2$).

On the other hand, a case by case analysis shows that none of the other BQDs,
classified in \cite[Section~10]{qsl3}, can give rise to an elliptic ATV-curve.
This proves Theorem~\ref{ASTtheo}.

\begin{rem}
The fact that, for the BQD $\L$ of Case~I.h, the ATV-curve of $\B_\L$ cannot
be elliptic was already visible in \eqref{ellnotell}, which was
obtained via an elimination procedure from Conditions~(\ref{coh}fg). In turn,
the latter were shown in \cite[Section~3]{qsl3} to be related to the
existence of an endomorphism of $V\o V$ satisfying the braid relation.

One would of course like to see a more direct and natural link, for an
arbitrary BQD $\L$, between the braid relation and the fact that the
ATV-curve of $\B_\L$ cannot be elliptic.
\end{rem}

\section*{Appendix: Multitwists preserve the Koszul property}
\renewcommand{\thesection}{A}
\setcounter{df}{0}

Let $\Ga$ be a monoid and $A=\bigoplus_{\ga\in\Ga}A_\ga$ a
$\Ga$-graded algebra.

Assume that to each $\ga\in\Ga$, we associate a graded automorphism
$\tau_\ga$ of $A$, in such a way that $\tau_{\ga\ga'}=\tau_\ga\tau_{\ga'}$
for all $\ga,\ga'\in\Ga$. (This is really only a special case of the notion
introduced in \cite{Zh}, but it will be sufficient for our purposes: in
Section~\ref{twistshapesec}, we take $\tau_{(k,l)}:=\tau_1^k\tau_2^l$ for
each $(k,l)\in\N^2$.)

Recall \cite{Zh} that the \emph{twisted} algebra $A^\tau$ is defined to be
the vector space $\bigoplus_{\ga\in\Ga}A_\ga$, endowed with the following new
multiplication:
\[
x*y:=x\tau_\ga(y)\qquad\text{for all $x\in A_\ga$, $y\in A$.}
\]
\begin{prop}
\label{Koszultwist}
Consider $A$ and $A^\tau$ as $\N$-graded algebras via some morphism
$h:\Ga\to\N$, and assume that they are generated by
$A_1(=A_1^\tau)=\bigoplus_{h(\ga)=1}A_\ga$. Then $A^\tau$ is Koszul if and
only if $A$ is Koszul.
\end{prop}

\begin{rem}
If $\Ga=\N$, then the categories of $\N$-graded modules of $A$ and of
$A^\tau$ are equivalent thanks to \cite[Theorem~3.1]{Zh}, and the Koszul
property is a homological property in this category, so the result is
immediate.

However, for arbitrary $\Ga$, \cite[Theorem~3.1]{Zh} shows that the
categories of $\Ga$-graded modules are equivalent, but this does not imply
that those of $\N$-graded modules are. Therefore, a proof is still needed.
\end{rem}

\begin{proof}[Proof of Proposition~\ref{Koszultwist}]
Of course, since $A=(A^\tau)^{\tau^{-1}}$, we only need to show one way, so
assume that $A$ is Koszul. Then $A$ is quadratic, say $A=T(A_1)/(R)$ with
$R\subset A_1\o A_1$. Moreover, for each $k\ge2$, the sublattice (w.r.t.\
$\cap$ and $+$) of $A_1^{\o k}$ generated by $R$ (that is, by all the
$A_1^{\o(i-1)}\o R\o A_1^{\o(k-i-1)}$) is distributive thanks to Backelin's
criterion \cite{Bath} (see also \cite[Lemma~4.5.1]{BGS}).

Now define $v:A_1^{\o k}\to A_1^{\o k}$ as follows: view $A_1^{\o k}$ as
the direct sum of all subspaces $A_{\ga_1}\o\dots\o A_{\ga_k}$ with
$h(\ga_1)=\dots=h(\ga_k)=1$, and set
\[
{v_{\mid}}_{A_{\ga_1}\o\dots\o A_{\ga_k}}
:=1_{A_{\ga_1}}\o\tau_{\ga_1}\o(\tau_{\ga_2}\tau_{\ga_1})
\o\dots\o(\tau_{\ga_{k-1}}\dots\tau_{\ga_1}).
\]
By construction, $A^\tau=T(A_1)/\bigl(v^{-1}(R)\bigr)$, so $A^\tau$ is again
quadratic. Moreover, the sublattice of $A_1^{\o k}$ generated by $v^{-1}(R)$
is the image under $v^{-1}$ of that generated by $R$, so it is still
distributive. Applying Backelin's criterion in the reverse direction, we
conclude that $A^\tau$ is Koszul.
\end{proof}

\end{document}